\documentclass[11pt]{article}
\usepackage{latexsym}
\usepackage{amsmath}
\usepackage{amsthm}
\usepackage{amssymb}

\newtheorem{theorem}{Theorem}[section]
\newtheorem{corollary}[theorem]{Corollary}
\newtheorem{lemma}[theorem]{Lemma}

\newtheorem{example}[theorem]{Example}

\DeclareMathOperator{\co}{co}

\newcommand{\covered}{<\hspace{-4pt}\cdot\hspace{2pt}}
\newcommand{\covers}{\hspace{2pt}\cdot\hspace{-4pt}>}
\DeclareMathOperator{\DP}{DP}
\DeclareMathOperator{\Kappa}{K}
\newcommand{\uuu}{{\bf u}}
\newcommand{\vvv}{{\bf v}}

\newcommand{\xxx}{{\bf x}}
\newcommand{\yyy}{{\bf y}}
\newcommand{\zzz}{{\bf z}}

\makeatletter
\@addtoreset{equation}{section}
\makeatother

\begin{document}

\title{The $f$-vector of the descent polytope}
\author{Denis Chebikin and Richard Ehrenborg}

\date{}

\maketitle

\begin{abstract}
For a positive integer $n$ and a subset
$S \subseteq [n-1]$,
the descent polytope~$\DP_{S}$ is the set of points 
$(x_{1}, \ldots, x_{n})$
in the
$n$-dimensional unit cube $[0,1]^{n}$
such that
$x_{i} \geq x_{i+1}$ if $i \in S$
and
$x_{i} \leq x_{i+1}$ otherwise.
First, we express the $f$-vector
as a sum over all subsets of $[n-1]$.
Second, we use certain factorizations
of the associated word 
over a two-letter alphabet to describe the $f$-vector.
We show that the $f$-vector is maximized when the
set $S$ is the alternating set $\{1,3,5, \ldots\} \cap [n-1]$.
We derive a generating function for $F_{S}(t)$,
written as a formal power series in two non-commuting variables
with coefficients in $\mathbb{Z}[t]$.
We also obtain the generating function for the 
Ehrhart polynomials of the descent polytopes.
\end{abstract}

\section{Introduction}
\label{section_introduction}

A classic topic in combinatorics is the study
of descent set statistics of permutations.
For a subset $S$ of the set
$[n-1] = \{1, \ldots, n-1\}$, the statistic $\beta(S)$ denotes the
number of permutations in the symmetric group
${\mathfrak S}_{n}$ with descent set $S$.
One well-known result is that 
the descent set statistic $\beta(S)$ is
maximized on the alternating sets:
$\{1,3, \ldots\} \cap [n-1]$
and
$\{2,4, \ldots\} \cap [n-1]$;
see~\cite{de_Bruijn, Niven, Readdy, Sagan_Yeh_Ziegler, Viennot}.

In this paper we study a class of polytopes $\DP_{S}$
which we call {\em descent polytopes}. They are indexed
by subsets $S$ of $[n-1]$, and the polytope corresponding to $S$
is the closure of the set of points of the unit hypercube $[0,1]^n$
whose coordinates, viewed as a sequence of $n$ numbers, have descent
set $S$. In general, these polytopes are
not simplicial nor simple: the polytope for $n=3$ and
$S = \{1\}$ is the Egyptian pyramid, that is, the square based
pyramid.
We show how to compute the $f$-vector of the descent polytope $\DP_S$.

Some invariants of descent polytopes are directly related to the
descent set statistic, and others exhibit analogous behavior.
For example, the volume of the descent polytope $\DP_{S}$ is given by the
descent set statistic $\beta(S)/n!$.  We also show in the paper that
the $f$-vector of the
descent polytope is entrywise maximized on the alternating set.
In order to prove this result we show that the entries of the $f$-vector
obey certain inequalities analogous to those satisfied by descent set
statistics;
see Theorem~\ref{theorem_inequality}.

One way to encode a subset $S$ of $[n-1]$
is by a word $\vvv_S$ in two letters, in our case $\xxx$ and $\yyy$.
Since $\vvv_S$ can be viewed as a non-commutative
monomial in the two variables $\xxx$ and $\yyy$,
this encoding suggests that one should work with non-commutative
generating functions. That is, to study a polytope invariant
$\phi$ of descent polytopes, one has to determine the
generating function
\begin{equation}
     \sum_{n \geq 1}
     \sum_{S \subseteq [n-1]}
      \phi(\DP_{S}) \cdot \vvv_{S}
  =
     \sum_{\vvv}
       \phi(\DP_{\vvv}) \cdot \vvv , 
\label{equation_phi}
\end{equation}
where we tacitly allow 
the descent polytopes
to be indexed by monomials.

In the case of the $f$-polynomial,
an encoding of the
$f$-vector, the generating function
in~\eqref{equation_phi} is a rational generating function;
see Theorem~\ref{theorem_Phi_formula}.
Furthermore, by expanding this rational function
we obtain a more concise expression
for the $f$-polynomial of the descent polytope~$\DP_{S}$.
This expression is in terms of a particular
type of factorizations of the monomial $\vvv_S$;
see Corollary~\ref{corollary_factorization}.

Descent polytopes are also lattice polytopes and hence
their Ehrhart polynomials are also of interest.
We also determine the non-commutative
generating function for Ehrhart polynomials of the descent polytopes;
see Theorem~\ref{theorem_Ehrhart}.
We note that even this power series is rational.

We end the paper with a few open questions and
directions for further research.

\section{An expression for the $f$-polynomial $F_{\vvv}$}
\label{section_F_S_summation}

For a set 
$S \subseteq [n-1] = \{1,2, \ldots, n-1\}$,
define the \emph{descent polytope}
$\DP_{S}$ to be the set of points
$(x_{1},\ldots,x_{n})$ in $\mathbb{R}^n$ such that
$0\leq x_{i} \leq 1$, and
$$
\left\{
\begin{array}{ll}
x_{i} \geq x_{i+1} & \mbox{if } i\in S, \\
x_{i} \leq x_{i+1} & \mbox{if } i\notin S.
\end{array}
\right.
$$
Thus $\DP_{S}$ is the \emph{order polytope} of
the ribbon poset $Z_{S} = \{z_{1},z_{2},\ldots,z_{n}\}$
defined by the cover relations $z_{i}\covers z_{i+1}$ if $i\in S$ and
$z_{i} \covered z_{i+1}$ if $i\notin S$; see~\cite{Stanley}.
It is clear that the set~$S$ and its complement
$\overline{S} = [n-1] - S$ yields the same descent polytope
up to an affine transformation. Also the
reverse set
$S^{rev} = \{n-i \: : \: i \in S\}$
give the same polytope.

Let $\xxx$ and $\yyy$ be two non-commuting variables.
For $S\subseteq [m]$, define $\vvv_{S} = \vvv_{1}\vvv_{2}\cdots\vvv_{m}$ where
$$
\vvv_{i} = \left\{\begin{array}{ll} \xxx & \mbox{if $i\notin S$,}\\
                                  \yyy & \mbox{if $i\in S$.}\end{array}\right.
$$
Since pairs 
$(n,S)$, where $S \subseteq[n-1])$, are in
bijective correspondence with $\xxx\yyy$-words
via $(n,S) \mapsto \vvv_{S}$,
it is natural to parameterize 
the descent polytopes and
their $f$-polynomials
by $\xxx\yyy$-words.
That is, we write $\DP_{\vvv}$ and $F_{\vvv}$,
where $\vvv = \vvv_{S}$ for some $S\subseteq [|\vvv|]$,
and $|\vvv|$ denotes the length of the word $\vvv$.
This notation has the advantage, as the $\xxx\yyy$-word
not only encodes the subset but also the dimension $n$.

For an $\xxx\yyy$-word $\vvv = \vvv_{1}\vvv_{2}\cdots\vvv_{n-1}$,
define the statistic $\kappa(\vvv)$ by 
$\kappa(\vvv) = 2 + \left|\{i\ :\ \vvv_{i} \neq \vvv_{i+1}\}\right|$
for $\vvv \neq 1$, and $\kappa(1) = 1$.
A direct observation is that
the number of facets of the descent polytope
is described by $\kappa$.
\begin{lemma}
The number of $(n-1)$-dimensional faces of the $n$-dimensional
descent polytope~$\DP_{\vvv}$ is given by
$$   f_{n-1}(\DP_{\vvv}) =  n-1 + \kappa(\vvv)  . $$
\end{lemma}
\begin{proof}
There are $n-1$ supporting hyperplanes of the form $x_{i} = x_{i+1}$
that each intersect the polytope in a facet.
The hyperplane $x_{i} = 1$ intersects the polytope in a facet if
one of the following three cases holds:
$\vvv_{i-1} \vvv_{i} = \xxx\yyy$;
$i=1$ and $\vvv_{1} = \yyy$; or
$i=n$ and $\vvv_{n} = \xxx$.
A similar statement holds for the hyperplane $x_{i} = 0$.
The lemma follows by adding these three statements.
\end{proof}

For an $\xxx\yyy$-word $\vvv = \vvv_{1}\vvv_{2}\cdots\vvv_{n-1}$
and a subset $T$ of $[n-1]$, define
$\vvv^{T}$  to be the subword
$\vvv^{T} = \vvv_{j_{1}} \vvv_{j_{2}} \cdots \vvv_{j_{k}}$,
where $T = \{j_{1} < j_{2} < \cdots < j_{k}\}$.
The following theorem provides a way to compute
the $f$-polynomial $F_{\vvv}$.

\begin{theorem}
Let $\vvv$ be an $\xxx\yyy$-word of length $n-1$.
Then the $f$-polynomial of the
descent polytope $\DP_{\vvv}$ is given by
$$
F_{\vvv}
      =
   1 
      +
   \sum_{T\subseteq [n-1]} 
           \left(\frac{t+1}{t}\right)^{\kappa(\vvv^{T})} 
         \cdot
           t^{|T|+1} .
$$
\label{theorem_F_formula}
\end{theorem}

\begin{proof}
For a face $\mathcal{F}$ of a polytope, let $\mathcal{F}^I$ denote the 
relative interior
of $\mathcal{F}$. Then the polytope is the disjoint union
of $\mathcal{F}^I$ taken over all faces $\mathcal{F}$, including the polytope
itself.

Recall that
the descent polytope $\DP_{\vvv}$ consists of all
points $(x_{1},\ldots,x_{n})\in\mathbb{R}^n$ belonging simultaneously to
the half spaces $x_{i} \geq 0$, $x_{i} \leq 1$ ($1\leq i \leq n$),
$x_{i} \leq x_{i+1}$ ($\vvv_{i} = \xxx$), and $x_{i} \geq x_{i+1}$
($\vvv_{i} = \yyy$).
A face $\mathcal{F}$ of $\DP_{\vvv}$ can be uniquely identified by
specifying which of these half spaces contain $\mathcal{F}$ on their
boundary hyperplanes, as long as the intersection of the whole
polytope and the specified boundary hyperplanes is non-empty.  Forming
the specification just for the half spaces of the form $x_{i}\leq
x_{i+1}$ or $x_{i}\geq x_{i+1}$ restricts the location of
$\mathcal{F}^I$ in $\mathbb{R}^n$ to the region defined by the
relations
\begin{equation}
x_{1} = x_{2} = \cdots = x_{j_{1}} \lessgtr
x_{j_{1}+1} = x_{j_{1}+2} = \cdots = x_{j_{2}} \lessgtr \cdots \lessgtr
x_{j_{k}+1} = x_{j_{k}+2} = \cdots = x_{n}
\label{equation_RT}
\end{equation}
for some $T = \{j_{1} < j_{2} <\cdots < j_{k}\}\subseteq [n-1]$, where the symbol
$\lessgtr$ denotes strict inequality: $x_{j_{i}} <x_{j_{i}+1}$
if $\vvv_{j_{i}} = \xxx$,
or $x_{j_{i}} > x_{j_{i}+1}$ if $\vvv_{j_{i}} = \yyy$.
Then $T$ is the set of indexes $j$ for which $\mathcal{F}$ does \emph{not} lie
entirely on the boundary hyperplane $x_j = x_{j+1}$ and thus 
the relative interior $\mathcal{F}^I$
is contained in the interior of the corresponding half space.
Let $\mathcal{R}(T)$ denote the intersection of the region defined by
(\ref{equation_RT}) and the hypercube $[0,1]^n$.
Each point $(x_{1},\ldots,x_{n})$ of $\DP_{\vvv}$
belongs to exactly one such region
$\mathcal{R}(T)$, namely, the one for $T = \{j\ |\ x_j\neq x_{j+1}\}$.
Thus we have the disjoint union
$$
\DP_{\vvv} = \bigsqcup_{T\subseteq [n-1]} \mathcal{R}(T).
$$
Let us show that the term
corresponding to $T\neq \varnothing$
in the expression in the statement of the theorem is
the contribution to $F_{\vvv}$ of the faces $\mathcal{F}$
of $\DP_{\vvv}$ for which $\mathcal{F}^I$ is contained in the region
$\mathcal{R}(T)$. In other words, we claim that for $T\neq \varnothing$ we
have
\begin{equation}
\sum_{\mathcal{F}\ :\ \mathcal{F}^I\subseteq\mathcal{R}(T)}
  t^{\dim \mathcal{F}}
=
  \left(\frac{t+1}{t}\right)^{\kappa(\vvv^{T})} 
\cdot
  t^{|T|+1}.
\label{equation_claim}
\end{equation}

Fix $\varnothing \neq T \subseteq [n-1]$. To select a particular face 
$\mathcal{F}$ from the
set of all faces with the property $\mathcal{F}^I \subseteq \mathcal{R}(T)$, we
need to complete the specification started above, that is, we must specify
which of the hyperplanes $x_{i} = 0,1$ contain $\mathcal{F}$, and we must make
sure that the intersection of the set of the specified hyperplanes and 
$\mathcal{R}(T)$ is non-empty. In terms of defining relations~(\ref{equation_RT}),
this task is equivalent to setting the common value of some of the ``blocks''
of coordinates
$(x_{1},\ldots,x_{j_{1}})$, $(x_{j_{1}+1}, \ldots, x_{j_{2}})$, \dots,
$(x_{j_{k}+1},\ldots,x_{n})$ to $0$ or $1$. Since the relations must remain
satisfiable by at least one point in $[0,1]^n$, only the blocks preceded
in~(\ref{equation_RT}) by $>$ (or nothing) and succeeded by $<$ (or nothing) can be set
to $0$. Similarly, only the blocks preceded by $<$ (or nothing) and succeeded
by $>$ (or nothing) can be set to~$1$. Thus each block can be set to at most
one of $0$ and $1$. The letters of the $\xxx\yyy$-word
$\vvv^{T} = \vvv_{j_{1}}\cdots \vvv_{j_{k}}$ encode
the inequality signs in~(\ref{equation_RT}) ($\xxx$ stands for $<$, and $\yyy$ stands
for $>$), so the number of blocks that can be set to $0$ or $1$ is the total
number of occurrences of $\xxx$ followed by $\yyy$, or $\yyy$ followed by
$\xxx$, in~$\vvv^{T}$, plus~$2$, as we also need to count the first and the
last blocks. In other words, the number of such blocks is $\kappa(\vvv^{T})$.

Observe that the dimension of the face of $\DP_{\vvv}$ obtained by
this specification procedure equals the number of blocks that have not been
set to $0$ or $1$: the common values of the coordinates in those blocks form
the ``degrees of freedom'' that constitute the dimension. Let us call such
blocks \emph{free}. The number of faces $\mathcal{F}$ with $\mathcal{F}^I
\subseteq \mathcal{R}(T)$ for which the specification procedure results in $m$
free blocks is
$$
    \binom{\kappa(\vvv^{T})}{|T|+1-m} ,
$$
the number of ways to choose $|T|+1-m$ blocks that are \emph{not} free out
of $\kappa(\vvv^{T})$ possibilities. Hence we have
\begin{eqnarray*}
\sum_{\mathcal{F}\ :\ \mathcal{F}^I\subseteq\mathcal{R}(T)}t^{\dim \mathcal{F}}
  & = &
\sum_{m=|T|+1-\kappa(\vvv^{T})}^{|T|+1}
     \binom{\kappa(\vvv^{T})}{|T|+1-m} \cdot t^{m}  \\
  & = &
t^{|T|+1-\kappa(\vvv^{T})} \cdot \sum_{\ell=0}^{\kappa(\vvv^{T})}
       \binom{\kappa(\vvv^{T})}{\ell} \cdot t^{\ell}  \\
  & = &
t^{|T|+1-\kappa(\vvv^{T})}\cdot (t+1)^{\kappa(\vvv^{T})},
\end{eqnarray*}
proving~(\ref{equation_claim}).

Finally, for $T=\varnothing$, we have
$\mathcal{R}(T) = \{0\leq x_{1}=\cdots=x_{n} \leq 1\}$, which is just the line
segment joining the two vertices $(0,\ldots,0)$ and $(1,\ldots,1)$ of
$\DP_{\vvv}$.
Thus the contribution of
$\mathcal{R}(T)$ to $F_{\vvv}$ is
$$
t+2 = 1
    + \left(\frac{t+1}{t}\right)^{\kappa(\vvv^\varnothing)} \cdot t.
$$
Adding this equation to the sum of~(\ref{equation_claim}) taken over the non-empty $T$
proves the theorem.
\end{proof}

Theorem~\ref{theorem_F_formula} yields
a combinatorial interpretation of the number of
vertices of the polytope $\DP_{\vvv}$. Call an $\xxx\yyy$-word $\vvv
= \vvv_{1} \vvv_{2}\cdots \vvv_{k}$ \emph{alternating}
if $\vvv_{i} \neq \vvv_{i+1}$ for all $1\leq i \leq k-1$.
Then we have the following corollary.

\begin{corollary}
For $\vvv$ an $\xxx\yyy$-word of length $n-1$,
the number of vertices of
the descent polytope $\DP_{\vvv}$ is one greater than
the number of subsets $T\subseteq [n-1]$ for which the word $\vvv^{T}$ is
alternating.
\label{corollary_number_vertices}
\end{corollary}

\begin{proof}
The number of vertices of $\DP_{\vvv}$ is the constant term of $F_{\vvv}$.
For the summand
corresponding to a subset $T\subseteq [n-1]$
in the formula of Theorem~\ref{theorem_F_formula},
the constant term is either $0$
or $1$, the latter being the case if and only if
$|T| + 1 - \kappa(\vvv^{T}) = 0$.
This condition is equivalent to $\vvv^{T}$ being alternating,
proving the corollary.
\end{proof}

As we mention in the introduction,
the descent set statistic $\beta(S)$
is maximized when $S$ is the alternating set.
The most elegant proof of this fact
uses the ${\bf c}{\bf d}$-index
of the simplex; see~\cite{Readdy}.
For an $\xxx\yyy$-word $\vvv$, let $\overline{\vvv}$ denote
the word obtained from $\vvv$ by replacing $\xxx$'s with~$\yyy$'s
and vice versa.
Then the following inequality holds:
\begin{equation}
   \beta(\uuu\yyy\xxx\vvv) > \beta(\uuu\yyy\yyy\overline{\vvv}) , 
\label{equation_classical_inequality}
\end{equation}
where we use $\xxx\yyy$-words to encode the sets.
In each of the
proofs~\cite{de_Bruijn, Niven, Sagan_Yeh_Ziegler, Viennot}
that the alternating
word maximizes the descent set,
the arguments rely on proving
the inequality~(\ref{equation_classical_inequality}).
However, the ${\bf c}{\bf d}$-index proof
gives a quick way to verify this inequality.
We now state a similar inequality for the $f$-vectors of descent
polytopes.

\begin{table}
$$
\begin{array}{c|c|c|c|c|c}
\uuu^{T} &
\vvv^{U} &
\kappa(\uuu^{T}\overline{\vvv}^{U}) &
\frac{Q_{T,U}(t)}{(t+1)^k\ t^{|T|+|U|+1-k}} &
\frac{\overline{Q}_{T,U}(t)}{(t+1)^{k}\ t^{|T|+|U|+1-k}} &
\frac{Q_{T,U}(t)-\overline{Q}_{T,U}(t)}{(t+1)^{k-1}\ t^{|T|+|U|+1-k}}\\[6pt]
\hline
&&&&&\\
\cdots\xxx & \xxx\cdots & k+1 & 1+t+(t+1)^{2} t^{-1} & (t+1)\ t^{-1}    &  (t+1)^{2}    \\[2pt]
           &            &     &    +\ (t+1)^{2}      & \cdot (1+2t+t^{2}) &             \\[9pt]
\cdots\xxx & \yyy\cdots & k-1 & 1+t+t              & (t+1)^{-1}\ t    &             \\[2pt]
           &            &     & +\ (t+1)^{2}       & \cdot \bigl(1+2(t+1)^{2}\ t^{-1}& t^{2} \\[2pt]
           &            &     &                    &  + (t+1)^{2}\bigr) &             \\[9pt]
\cdots\yyy & \xxx\cdots & k-1 & 1+t+t+t^{2}          & (t+1)^{-1}\ t    &  (t+1)^{2}    \\[2pt]
           &            &     &                    & \cdot(1+2t+t^{2})  &             \\[9pt]
\cdots\yyy & \yyy\cdots & k+1 & 1+(t+1)^{2}\ t^{-1}  & (t+1)\ t^{-1}    &  (t+1)^{2}    \\[2pt]
           &            &     &  +\ t + (t+1)^{2}    & \cdot (1+2t+t^{2}) &             \\[9pt]
1          & \xxx\cdots &  k  & 1+t+(t+1)          & 1 + 2t + t^{2}     &  (t+1)^{2}    \\[2pt]
           &            &     & +\ (t+1)\ t        &                  &             \\[9pt]
1          & \yyy\cdots &  k  & 1+(t+1)            & 1 + 2(t+1)       &  t^{2}+t      \\[2pt]
           &            &     & +\ t + (t+1)^{2}     & +\ (t+1)\ t      &             \\[9pt]
\cdots\xxx & 1          &  k  & 1+t+(t+1)          & 1 + 2(t+1)       &  t^{2}+t      \\[2pt]
           &            &     & +\ (t+1)^{2}         & +\ (t+1)\ t      &             \\[9pt]
\cdots\yyy & 1          &  k  & 1+(t+1)            & 1 + 2t + t^{2}     &  (t+1)^{2}    \\[2pt]
           &            &     & +\ t + (t+1)\ t    &                  &             \\[9pt]
1          & 1          & k=1 & 1+(t+1)            & 1+ 2(t+1)        &  (t+1)^{2}    \\[2pt]
           &            &     & +\ (t+1) + (t+1)^{2} & +\ (t+1)\ t      &
\end{array}
$$
\caption{Calculations for the proof of Theorem~\ref{theorem_inequality}.}
\label{table_one}
\end{table}

\begin{theorem}
Let $\uuu$ and $\vvv$ be two $\xxx\yyy$-words
such that the sum of their lengths is $n-3$,
that is, $|\uuu| + |\vvv| = n-3$.
Then the difference
\begin{equation}
F_{\uuu\yyy\xxx\vvv}(t) - F_{\uuu\yyy\yyy\overline{\vvv}}(t)
\label{equation_difference}
\end{equation}
has positive coefficients at $1, t, \ldots, t^{n-1}$.
That is, for $0 \leq i \leq n-1$
the descent polytope
$\DP_{\uuu\yyy\xxx\vvv}$
has more faces of dimension $i$ than
the descent polytope $\DP_{\uuu\yyy\yyy\overline{\vvv}}$.
\label{theorem_inequality}
\end{theorem}
\begin{proof}
Let $|\uuu|=m$ and $|\vvv|=n-m-3$. For $T\subseteq [m]$ and
$U\subseteq [n-m-3]$, define
\begin{equation}
Q_{T,U}(t) 
   :=
\sum_{E\subseteq \{1,2\}} 
    \left(\frac{t+1}{t}\right)^{\kappa(\uuu^{T}(\yyy\xxx)^E\vvv^{U})}
  \cdot
    t^{|T|+|U|+|E|+1} .
\label{equation_Q_T_U_definition}
\end{equation}
Thus $Q_{T,U}(t)$ is the sum of four of the terms in the summation formula
for $F_{\uuu\yyy\xxx\vvv}(t)$ given by Theorem~\ref{theorem_F_formula},
corresponding
to fixed choices of letters drawn from $\uuu$ and from $\vvv$.
Similarly, let us define
$$
\overline{Q}_{T,U}(t) 
  :=
\sum_{E\subseteq \{1,2\}} 
    \left(\frac{t+1}{t}\right)^{\kappa(\uuu^{T}(\yyy\yyy)^E\overline{\vvv}^{U})}
  \cdot
    t^{|T|+|U|+|E|+1}.
$$ 
Note that $Q_{T,U}(t)$ depends only on $|T|$, $|U|$,
$\kappa(\uuu^{T}\vvv^{U})$, the last letter of~$\uuu^{T}$,
and the first letter
of~$\vvv^{U}$, and not on the
particular choice of the
remaining letters of
$\uuu^{T}$ and $\vvv^{U}$. Thus to
show that the difference $Q_{T,U}(t)-\overline{Q}_{T,U}(t)$ is
a polynomial with non-negative coefficients,
it suffices to consider $9$ cases
corresponding to $\uuu^{T}$ (respectively, $\vvv^{U}$) ending
(respectively, beginning) with $\xxx$ or $\yyy$,
or being equal to the empty word $1$.

We summarize our calculations in Table~\ref{table_one}.
We denote $\kappa(\uuu^{T}\vvv^{U})$ by
$k$, and we divide each polynomial by the common factor
$(t+1)^k \cdot t^{|T|+|U|+1-k}$.
In the fourth column, which corresponds to $Q_{T,U}$, the four summands
represent the results of inserting~$1$, $\xxx$, $\yyy$, and $\yyy\xxx$ between
$\uuu^{T}$ and~$\vvv^{U}$. For example, if~$\uuu^{T}$
ends with an $\xxx$ and
$\vvv^{U}$ begins with an $\xxx$, then inserting $\yyy$ increases the value
of the statistic $\kappa$ by~$2$, thus contributing a factor of
$$
     \left( \frac{t+1}{t} \right)^{2} \cdot t
   = 
     (t+1)^{2} \cdot t^{-1}
$$
to the corresponding term of~(\ref{equation_Q_T_U_definition}).
Similarly, the entries in the fifth column
consist of a factor resulting from a different value of $\kappa$ for the word
$\uuu^{T}\overline{\vvv}^{U}$ times the contributions of inserting $1$, $\yyy$
(counted twice), and $\yyy^{2}$ between $\uuu^{T}$ and $\overline{\vvv}^{U}$.

We conclude that in every case
the quotient
$\frac{Q_{T,U}(t)-\overline{Q}_{T,U}(t)}
      {(t+1)^{k-1} \cdot t^{|T|+|U|+1-k}}$
is a polynomial of degree~$2$
with non-negative coefficients.
Hence the difference $Q_{T,U}(t)-\overline{Q}_{T,U}(t)$
is a polynomial with non-negative coefficients.
Summing over all possible pairs
$(T,U)$ yields that the difference
in equation~(\ref{equation_difference})
has non-negative coefficients.
More specifically, the polynomial
$Q_{T,U}(t)-\overline{Q}_{T,U}(t)$
has degree
$(k-1) + (|T|+|U|+1-k) + 2 = |T|+|U|+2$.
This degree can attain
any integer value between $2$ and $n-1$.
Thus the leading terms of these
differences contribute positively to the coefficients
of $t^{2}, t^{3}, \ldots, t^{n-1}$
in the 
difference~(\ref{equation_difference}).
Furthermore,
in the case $T=U=\varnothing$ we have
$Q_{T,U}(t)-\overline{Q}_{T,U}(t) = (t+1)^{2}$,
which yields a positive contribution to the constant and the
linear terms of the overall difference. The proof is now complete.
\end{proof}

Let~$\zzz_{n}$ be the alternating word of length $n$
starting with the letter $\xxx$.
Then $\overline{\zzz_{n}}$ is the alternating word
beginning with $\yyy$. That is, the two alternating words are
$$    \zzz_{n}
    =
      \underbrace{\xxx\yyy\xxx\cdots}_{n}
                  \:\:\:\: \mbox{ and } \:\:\:\:
      \overline{\zzz_{n}}
    =
      \underbrace{\yyy\xxx\yyy\cdots}_{n} .
$$
We now have the maximization result for the $f$-vector
of descent polytopes.
\begin{corollary}
The $f$-vector of the two descent polytopes
$\DP_{\zzz_{n-1}}$ and $\DP_{\overline{\zzz_{n-1}}}$
is maximal among the $f$-vectors of all descent polytopes
of dimension $n$. That is, for each $0 \leq i \leq n-1$,
the polytope
$\DP_{\zzz_{n-1}}$ has more faces of 
dimension $i$ than the descent polytope $\DP_{\vvv}$ of dimension $n$
for a non-alternating word $\vvv$.
\label{corollary_alternating_pattern_max}
\end{corollary}

\section{The power series $\Phi(\xxx,\yyy)$}
\label{section_Phi_formula}

We now derive a non-commutative generating function
$\Phi(\xxx,\yyy)$ for the $f$-polynomial~$F_{\vvv}$,
which belongs to the ring 
$\Phi(\xxx,\yyy)\in \mathbb{Z}[t]\langle\langle \xxx,\yyy \rangle\rangle$.
We define the power series $\Phi(\xxx,\yyy)$ by 
$$
\Phi(\xxx,\yyy) = \sum_{\vvv} F_{\vvv} \cdot \vvv  ,  
$$ 
where the sum is over all $\xxx\yyy$-words $\vvv$.
Since we have the symmetry $F_{\vvv} = F_{\overline{\vvv}}$,
we obtain that $\Phi(\xxx,\yyy)$ is symmetric
with respect to $\xxx$ and $\yyy$, that is,
$$
\Phi(\xxx,\yyy) = \Phi(\yyy,\xxx) .
$$

Let $\vvv$ be an $\xxx\yyy$-word $\vvv_{1}\vvv_{2}\cdots \vvv_{n-1}$.
Consider the following polynomials:
\begin{eqnarray*}
K_{\vvv}(t)
  & := & 
\sum_{T\subseteq [n-1]\ :\ \vvv_{j_{1}}=\xxx}
      \left(\frac{t+1}{t}\right)^{\kappa(\vvv^{T})}
   \cdot
      t^{|T|+1}, \\
L_{\vvv}(t)
  & := &
\sum_{T\subseteq [n-1]\ :\ \vvv_{j_{1}}=\yyy}
      \left(\frac{t+1}{t}\right)^{\kappa(\vvv^{T})}
   \cdot
      t^{|T|+1},
\end{eqnarray*}
where $\vvv_{j_{1}}$ denotes the first letter of the word
$\vvv^{T} = \vvv_{j_{1}}\vvv_{j_{2}}\cdots \vvv_{j_{k}}$, as in the
notation of Theorem~\ref{theorem_F_formula}.
Since $\vvv^{T}$ begins with either 
$\xxx$ or $\yyy$ unless $T=\varnothing$, we have
\begin{equation}
F_{\vvv} = K_{\vvv} + L_{\vvv} + t + 2,
\label{equation_F_from_KL}
\end{equation}
where $t+2$ is the $f$-polynomial of $\DP_{1}$,
the line segment.
We continue with a lemma that relates 
the two polynomials $K_{\vvv}$ and $L_{\vvv}$.

\begin{lemma}
For an $\xxx\yyy$-word $\vvv$
the following four equalities hold:
$$
\begin{array}{l}
K_{\yyy\vvv} = K_{\vvv}, \\
L_{\xxx\vvv} = L_{\vvv}, \\
K_{\xxx\vvv} = L_{\yyy\vvv} = (t+1) \cdot (K_{\vvv} + L_{\vvv} + t + 1).
\end{array}
$$
\label{lemma_KL_recurrence}
\end{lemma}
\begin{proof}
For an integer $i$ and a set $U\subseteq \mathbb{Z}$, let $U+i$ denote the set
obtained by adding $i$ to each element of $U$.
Also let $\vvv = \vvv_{1}\vvv_{2}\cdots\vvv_{n-1}$,
where each $\vvv_{i}$ is either $\xxx$ or $\yyy$.

Clearly, $(\yyy\vvv)^{T}$ begins
with $\xxx$ if and only if $1\notin T$ and $\vvv^{T-1}$ begins with $\xxx$,
in which case $(\yyy\vvv)^{T} = \vvv^{T-1}$. Hence $K_{\yyy\vvv} = K_{\vvv}$.

Now, $(\xxx\vvv)^{T}$ begins with $\xxx$ if and only if either $1\in T$, or else
$1 \notin T$ and $\vvv^{T-1}$ begins with $\xxx$.
In the former case, we have $T = \{1 < j_{1} + 1 < j_{2} + 1 < \cdots < j_{k}+1\}$,
and $(\xxx\vvv)^{T} = \xxx \vvv_{j_{1}} \vvv_{j_{2}} \cdots \vvv_{j_{k}}$.
Set $U = (T - \{1\}) - 1 = \{j_{1} < \cdots <j_{k}\}$.
Then $\kappa\left((\xxx\vvv)^{T}\right) = \kappa\left(\vvv^{U}\right)$
if $\vvv_{j_{1}}=\xxx$, and 
$\kappa\left((\xxx\vvv)^{T}\right) = \kappa\left(\vvv^{U}\right) + 1$ if
$\vvv_{j_{1}}=\yyy$.
Hence
\begin{eqnarray}
\label{equation_K_a_v_1}
\sum_{1\in T\subseteq [n]}
\left(\frac{t+1}{t}\right)^{\kappa\left((\xxx\vvv)^{T}\right)}
 \cdot
t^{|T|+1} 
  & = &
(t+1)^{2} + 
t\cdot \sum_{U\ :\ \vvv_{j_{1}} = \xxx} 
\left(\frac{t+1}{t}\right)^{\kappa\left(\vvv^{U}\right)} 
\cdot t^{|U|+1} \\
  &   &
+ (t+1)\cdot \sum_{U\ :\ \vvv_{j_{1}}=\yyy} 
\left(\frac{t+1}{t}\right)^{\kappa\left(\vvv^{U}\right)} 
\cdot
t^{|U|+1}
\nonumber \\
  & = &
(t+1)^{2} + t \cdot K_{\vvv} + (t+1) \cdot L_{\vvv},
\nonumber
\end{eqnarray}
where the first term $(t+1)^{2}$ corresponds to $T = \{1\}$ and
$U=\varnothing$. In the case where $1\notin T$ and $\vvv^{T-1}$ begins with
$\xxx$ we have, as before, $(\xxx\vvv)^{T} = \vvv_{j_{1}} \vvv_{j_{2}}\cdots 
\vvv_{j_{k}} = \vvv^{T-1}$, and hence
\begin{equation}
\sum_{T\ :\ \vvv_{j_{1}}=\xxx}
\left(\frac{t+1}{t}\right)^{\kappa\left((\xxx\vvv)^{T}\right)}
\cdot
t^{|T|+1}
=
\sum_{\vvv_{j_{1}}=\xxx}
\left(\frac{t+1}{t}\right)^{\kappa\left(\vvv^{T-1}\right)}
\cdot
t^{|T-1|+1}
=
K_{\vvv}.
\label{equation_K_a_v_2}
\end{equation}

Adding~(\ref{equation_K_a_v_1}) and~(\ref{equation_K_a_v_2}) yields
$$
K_{\xxx\vvv} = (t+1) \cdot (K_{\vvv} + L_{\vvv} + t + 1).
$$
The relations for $L_{\xxx\vvv}$ and $L_{\yyy\vvv}$ follow from symmetry
that arises from exchanging the variables $\xxx$ and~$\yyy$.
\end{proof}

Starting with $K_{1} = L_{1}= 0$, one can use
Lemma~\ref{lemma_KL_recurrence}
to recursively compute $K_{\vvv}$ and $L_{\vvv}$, and
hence~$F_{\vvv}$, from~(\ref{equation_F_from_KL}).
Recall the generating power series
$$
\Phi(\xxx,\yyy) = \sum_{\vvv} F_{\vvv} \cdot \vvv,
$$
where the sum is over all $\xxx\yyy$-words, including the empty
word $\vvv = \vvv_{\varnothing} = 1$.
Define the two generating power series
\begin{eqnarray*}
\Kappa(\xxx,\yyy)
  & := &
\sum_{\vvv} K_{\vvv} \cdot \vvv , \\
\Lambda(\xxx,\yyy)
  & := &
\sum_{\vvv} L_{\vvv} \cdot \vvv .
\end{eqnarray*}
{}From the definitions of $K_{\vvv}$ and $L_{\vvv}$ it follows that
$K_{\vvv} = L_{\overline{\vvv}}$.
By the symmetry in the two variables $\xxx$ and $\yyy$
we have that
$$
\Lambda(\xxx,\yyy) = \Kappa(\yyy,\xxx) .
$$ 
Then, by~(\ref{equation_F_from_KL}), we have
\begin{eqnarray}
\nonumber
\Phi(\xxx,\yyy)
  & = &
\Kappa(\xxx,\yyy) + \Lambda(\xxx,\yyy)+(t+2) \cdot \sum_{\vvv} \vvv \\
\nonumber
  & = &
\Kappa(\xxx,\yyy) + \Kappa(\yyy,\xxx) 
+ (t+2) \cdot \sum_{r\geq 0} (\xxx+\yyy)^{r} \\
  & = &
 \Kappa(\xxx,\yyy) + \Kappa(\yyy,\xxx) + (t+2) \cdot \frac{1}{1-\xxx-\yyy} .
\label{equation_Phi_from_Kappa}
\end{eqnarray}
Using the equations in Lemma~\ref{lemma_KL_recurrence}
and recalling that $K_{1} = 0$  we obtain
\begin{eqnarray*}
\Kappa(\xxx,\yyy)
  & = &
\sum_{\vvv} K_{\xxx\vvv} \cdot \xxx\vvv 
  +
\sum_{\vvv} K_{\yyy\vvv} \cdot \yyy\vvv \\
  & = &
(t+1) \cdot \xxx \cdot \sum_{\vvv} (K_{\vvv}+L_{\vvv}+t+1) \cdot \vvv
  +
\yyy \cdot \sum_{\vvv} K_{\vvv} \cdot \vvv \\
  & = &
   (t+1) \cdot \xxx \cdot
    \left(\Kappa(\xxx,\yyy) + \Lambda(\xxx,\yyy) 
              + (t+1)\cdot \frac{1}{1-\xxx-\yyy} \right)
    + \yyy \cdot \Kappa(\xxx,\yyy)\\
  & = &
 (t+1) \cdot \xxx \cdot \left(\Phi(\xxx,\yyy) - \frac{1}{1-\xxx-\yyy}\right) +
\yyy \cdot \Kappa(\xxx,\yyy),
\end{eqnarray*}
where the last step is by equation~(\ref{equation_Phi_from_Kappa}).
Rearranging terms we have
$$
\Kappa(\xxx,\yyy)
  =
    (t+1) \cdot (1-\yyy)^{-1} \cdot \xxx 
  \cdot
    \left(\Phi(\xxx,\yyy) - \frac{1}{1 - \xxx-\yyy}\right) .
$$
Adding this equation and its symmetric version obtained
by exchanging $\xxx$ and $\yyy$ one has
\begin{eqnarray*}
\Kappa(\xxx,\yyy)+\Kappa(\yyy,\xxx)
  & = &
(t+1)
  \cdot 
\left((1-\yyy)^{-1} \cdot \xxx + (1-\xxx)^{-1} \cdot \yyy\right)
  \cdot
\left(\Phi(\xxx,\yyy) - \frac{1}{1-\xxx-\yyy}\right) ,
\end{eqnarray*}
using the symmetry $\Phi(\xxx,\yyy) = \Phi(\yyy,\xxx)$.
Now using equation~(\ref{equation_Phi_from_Kappa})
we can solve for $\Phi(\xxx,\yyy)$ and arrive at the following theorem.
\begin{theorem}
The generating power series $\Phi(\xxx,\yyy)$ is given by
$$
\Phi(\xxx,\yyy)
  =
\left(1 + \frac{t+1}{1 - (t+1) \cdot \left((1-\yyy)^{-1} \cdot \xxx + (1-\xxx)^{-1} \cdot \yyy\right)}
\right)\cdot
\frac{1}{1-\xxx-\yyy}.
$$
\label{theorem_Phi_formula}
\end{theorem}

\begin{corollary}
For an $\xxx\yyy$-word $\vvv$ the $f$-vector of the descent polytope
$\DP_{\vvv}$ is given by the sum
$$   F_{\vvv}(t)
  =
     1 + \sum_{(\uuu_{1}, \ldots, \uuu_{k-1}, \uuu_{k})} (t+1)^{k} ,
$$
where the sum ranges over all factorizations of the
word $\vvv = \uuu_{1} \cdots \uuu_{k-1} \cdot \uuu_{k}$
such that each of the factors
$\uuu_{1}, \ldots, \uuu_{k-1}$
are of the form $\xxx^{i} \yyy$ or $\yyy^{i} \xxx$,
where $i \geq 0$,
and there is no condition on the last factor~$\uuu_{k}$.
\label{corollary_factorization}
\end{corollary}
\begin{proof}
Rewrite Theorem~\ref{theorem_Phi_formula} as
\begin{eqnarray*}
      \Phi(\xxx,\yyy)
  & = &
     \frac{1}{1-\xxx-\yyy} 
   +
     \frac{1}{1 - (t+1) \cdot \left((1-\yyy)^{-1} \cdot \xxx + (1-\xxx)^{-1} \cdot \yyy\right)}
       \cdot
     \frac{t+1}{1-\xxx-\yyy} \\
  & = &
       \sum_{\vvv}
         \vvv
   +
     \sum_{j \geq 0}
       \left(
           (t+1)
         \cdot
           \sum_{i \geq 0} 
             \left(\yyy^{i} \xxx + \xxx^{i} \yyy\right)
       \right)^{j}
   \cdot
     (t+1)
   \cdot
       \sum_{\vvv}
         \vvv    ,
\end{eqnarray*}
where in both sums $\vvv$ ranges over all $\xxx\yyy$-words.
The corollary follows by reading the generating function.
\end{proof}

\begin{example}
{\rm
Consider the $5$-dimensional
descent polytope $\DP_{\vvv}$
where $\vvv = \xxx\yyy\yyy\xxx$.
We have the following list of $11$ factorizations:
$$ \begin{array}{l l l l l l l l l}
\vvv & = & \xxx\yyy\yyy\xxx \\
     & = & \xxx \cdot \yyy\yyy\xxx 
     & = & \xxx\yyy \cdot \yyy\xxx \\
     & = & \xxx \cdot \yyy \cdot \yyy\xxx 
     & = & \xxx \cdot \yyy\yyy\xxx \cdot 1 
     & = & \xxx\yyy \cdot \yyy \cdot \xxx 
     & = & \xxx\yyy \cdot \yyy\xxx \cdot 1 \\
     & = & \xxx \cdot \yyy \cdot \yyy \cdot \xxx 
     & = & \xxx \cdot \yyy \cdot \yyy\xxx \cdot 1 
     & = & \xxx\yyy \cdot \yyy \cdot \xxx \cdot 1 \\
     & = & \xxx \cdot \yyy \cdot \yyy \cdot \xxx \cdot 1
   \end{array} $$
Hence the $f$-polynomial of the polytope $\DP_{\xxx\yyy\yyy\xxx}$
is given by
\begin{eqnarray*}
     F_{\xxx\yyy\yyy\xxx}
  & = &
     1 + (t+1) + 2 \cdot (t+1)^{2}
               + 4 \cdot (t+1)^{3}
               + 3 \cdot (t+1)^{4}
               + (t+1)^{5}  \\
  & = &
     12 
   + 34 \cdot t 
   + 42 \cdot t^{2}
   + 26 \cdot t^{3}
   +  8 \cdot t^{4}
   +          t^{5}  .
\end{eqnarray*}
}
\end{example}

For the alternating word $\zzz_{n-1}$ we can say more
about the associated descent polytope.
The number of vertices of
$\DP_{\zzz_{n-1}}$ is the Fibonacci number $F_{n+2}$;
see for instance~\cite[Exercise~1.14e]{Stanley_EC_I}.
More generally, the $f$-vector of
$\DP_{\zzz_{n-1}}$ is given by the next result.
\begin{corollary}
The $f$-polynomial of the $n$-dimensional descent
polytope $\DP_{\zzz_{n-1}}$
is described by
$$     F_{\zzz_{n-1}}
    =
       1
     +
       \sum_{(c_{1}, c_{2}, \ldots, c_{k})}
                 (t+1)^{k}   ,   $$
where the sum is over all compositions of $n$ such that
all but the last part is less than or equal to $2$,
that is,
$c_{1}, \ldots, c_{k-1} \in \{1,2\}$.
\end{corollary}
\begin{proof}
The only factors of the alternating word $\zzz_{n-1}$
of the form $\xxx^{i} \yyy$ or $\yyy^{i} \xxx$
have $i = 0,1$. Hence it is enough
to record the length of each factor $\uuu_{i}$,
that is, $d_{i} = |\uuu_{i}|$. Thus we are summing
over vectors of non-negative integers
$(d_{1}, \ldots, d_{k})$
such that
the sum of the entries is $n-1$
and
$d_{1}, \ldots, d_{k-1} \in \{1,2\}$
and $d_{k} \geq 0$. By adding one to the last entry $d_{k}$
we have a composition of $n$.
\end{proof}

This corollary yields the generating function
\begin{equation}
      \sum_{n \geq 1}
           F_{\zzz_{n-1}} \cdot x^{n}
 =
      \frac{x}{1-x}
    +
        \frac{1}{1 - (t+1) \cdot (x + x^{2})}
      \cdot 
        (t+1)
      \cdot 
        \frac{x}{1-x}   .  
\label{equation_alternating}
\end{equation}
Setting $t=0$ in this generating function and adding constant $1$
yields $(1+x)/(1-x-x^{2})$,
the generating function for the Fibonacci numbers as expected.

\section{A generating power series for the Ehrhart polynomials
of $\DP_S$}
\label{section_Ehrhart}

Besides the $f$-vector,
another geometric invariant of a polytope
is the \emph{Ehrhart polynomial}.
As a function of a non-negative integer $r$,
the Ehrhart polynomial of a lattice polytope $P$
is the number of lattice points in the dilation
$r \cdot P$.
Ehrhart's fundamental result is that this function
is a polynomial in $r$.
In the case of the $n$-dimensional descent polytope $\DP_{S}$,
the Ehrhart polynomial $\iota_{S}(r)$
counts the number of lattice points satisfying
the inequalities
$0 \leq x_1, x_2, \ldots, x_n \leq r$,
$x_{i} \geq x_{i+1}$ for $i \in S$, and $x_{i} \leq x_{i+1}$
for $i \notin S$.
In this section we derive the generating power series
$$
I(r;\xxx,\yyy) := \sum_{\vvv} \iota_{\vvv}(r)\cdot \vvv,
$$
As before, we adopt the shorthand
$\iota_{\vvv_{S}} = \iota_{S}$.

Let us call an element $w$ of the set
$\{0,1,\ldots,r\}^n$ an \emph{$r$-word} of
length $n$.
Define the descent set $D(w)$ of $w=(w_{1},w_{2},\ldots,w_{n})$
to be the set of positions $i$ such that $w_{i}>w_{i+1}$.  For an
$\xxx\yyy$-word $\vvv$ of length~$n-1$, let $\beta(r,\vvv)$ be the number of
$r$-words $w$ of length $n$ such that $D(w)$ is encoded by $\vvv$ (that
is, $\vvv=\vvv_{D(w)}$).  Note that $\beta(r,\vvv)$ is not quite the
Ehrhart polynomial $\iota_\vvv$, as it only counts those integer points
of $r\DP_\vvv$ with strict descents.  Still, a
generating power series for $\beta(r,\vvv)$ is the first step in our
computation of~$I(r;\xxx,\yyy)$.

For an $\xxx\yyy$-word of length $n-1$, let $\alpha(r,\vvv)$ be the number
of $r$-words of length $n$ such that $D(w)$ is contained in the subset
of $[n-1]$ encoded by $\vvv$.  Fix $r\geq 0$, $n>0$, and $S\subseteq
[n-1]$.  Write the word~$\vvv_{S}$ as
$$
\vvv_{S}
  =
\xxx^{g_{1}-1} \yyy \xxx^{g_{2}-1} \yyy \cdots \yyy \xxx^{g_{k}-1},
$$
so that $g=(g_{1},g_{2},\ldots,g_{k})$ is the composition $\co(S)$
of $n$ associated to $S$.
To construct a word counted by $\alpha(r,\vvv)$,
one needs to choose, for every $i$, a multiset of elements of
$\{0,1,\ldots,r\}$ that go into the ``block''
corresponding to the part $g_{i}$ of $g$,
put them in (weakly) increasing order, and concatenate the blocks.
There are $\binom{r+g_{i}}{g_{i}}$ ways of choosing the elements for the
$i$-th block.  Thus if we define
$$
Q_{r}(x) := \sum_{j\geq 1} \binom{r+j}{j} \cdot x^{j-1}
          = x^{-1} \cdot \left((1-x)^{-r-1} - 1\right) 
$$
then the generating power series for $\alpha(r,\vvv)$ is
\begin{eqnarray}
A(r;\xxx,\yyy)
  & := &
\sum_\vvv \alpha(r,\vvv)\cdot \vvv 
\nonumber \\
  & = &
\sum_{k \geq 1}
    Q_{r}(\xxx)
  \cdot
\left(\yyy \cdot Q_{r}(\xxx)\right)^{k-1}
\nonumber \\
  & = &
    Q_{r}(\xxx)
  \cdot
\left(1 - \yyy \cdot Q_{r}(\xxx)\right)^{-1} ,
\label{A_rxy}
\end{eqnarray}
where $k$ runs through all possible numbers of parts of the
composition $g = \co(S)$, and the index $j$ in the definition of
$Q_{r}$ corresponds to the choice of each part size.

A standard application of the inclusion-exclusion principle yields
\begin{equation}
B(r;\xxx,\yyy) :=  \sum_\vvv \beta(r;\xxx,\yyy)\ \vvv
                =  A(r;\ \xxx-\yyy,\ \yyy).
\label{B_rxy}
\end{equation}

The final step is the following claim:

\begin{lemma}
The generating function of the Ehrhart polynomials
of the descent polytopes is expressed in terms
of $B(r;\xxx,\yyy)$ as
$$
I(r;\xxx,\yyy) = (1-\yyy)^{-1}\cdot B\left(r;\ \xxx(1-\yyy)^{-1},
                                      \ \yyy(1-\yyy)^{-1}\right).
$$
\end{lemma}
\begin{proof}
For an $\xxx\yyy$-word $\vvv = \vvv_{1} \vvv_{2} \cdots \vvv_{n-1}$ 
and an integer point
$p=(p_{1}, p_{2}, \ldots, p_{n})$ in $r\DP_\vvv$, define $\varphi(\vvv, p) = (\uuu, q)$,
where the $\xxx\yyy$-word $\uuu$ and the $r$-word $q$ are
obtained from $\vvv$ and $p$ as follows: for each index $i\in [n-1]$ such
that $p_{i} = p_{i+1}$ and $\vvv_{i} = \yyy$, remove the letter at the $i$-th
position from $\vvv$ as well as the $i$-th coordinate from $p$.  For example,
if $p = (2, 3, 1, 1, 1, 1, 1, 4)$ and $\vvv = \xxx\yyy\yyy\yyy\xxx\yyy\xxx$,
then the removal should be done for $i = 3, 4, 6$, so $q = (2, 3, 1, 1, 4)$
and $\uuu = \xxx\yyy\xxx\xxx$. 

Let $B(r,\uuu)$ be the set of $r$-words with descent set encoded by $\uuu$,
so that $\beta(r,\uuu) = \left|B(r,\uuu)\right|$.
Note that if
$\varphi(\vvv,p) = (\uuu,q)$
then $q \in B(r,\uuu)$.  For a
fixed $\uuu$ and $q\in B(r,\uuu)$, the inverse image
 $\varphi^{-1}(\uuu,q)$
can be obtained by performing the following operation in all possible ways:
start with $\uuu$, insert an arbitrary number of $\yyy$'s (maybe none)
in each of the gaps between consecutive letters of $\uuu$, before the
first letter of $\uuu$, and after the last letter of $\uuu$,
 and for each coordinate $q_{i}$ of $q$, insert as many
copies of $q_{i}$ before that coordinate as the number of $\yyy$'s that were
inserted before the $i$-th letter of $\uuu$
(for $i=n$, use the number of $\yyy$'s inserted after the last letter
of $\uuu$).  The resulting $\xxx\yyy$-word
$\vvv$ and integer point $p$ satisfy $\varphi(\vvv,p) = (\uuu,q)$.
In terms generating functions we have
\begin{eqnarray}
\sum_{(\vvv,p) \in \varphi^{-1}(\uuu,q)} \vvv
  & = &
(1+\yyy+\yyy^{2}+\cdots)
  \cdot
\uuu_{1}
  \cdot
(1+\yyy+\yyy^{2}+\cdots)
  \cdot
\uuu_{2}
  \cdots
\uuu_{m}
  \cdot
(1+\yyy+\yyy^{2}+\cdots)
\nonumber \\
  & = &
(1-\yyy)^{-1}
  \cdot
\uuu_{1}
  \cdot
(1-\yyy)^{-1}
  \cdot
\uuu_{2}
  \cdots
\uuu_{m}
  \cdot
(1-\yyy)^{-1} ,
\label{fixed_u_q}
\end{eqnarray}
where $\uuu = \uuu_{1}\uuu_{2}\cdots \uuu_{m}$ and $\uuu_{i}\in\{\xxx,\yyy\}$.
Consider the sum of equation (\ref{fixed_u_q})
over all pairs $(\uuu,q)$ such that $q\in B(r,\uuu)$.
The left-hand side
of the resulting identity is
$$
\sum_{\left\{(\vvv,p)\ :\ p \in r\DP_\vvv\right\}} \vvv 
= \sum_\vvv \iota_\vvv(r) \cdot \vvv = I(r;\xxx,\yyy).
$$
The right-hand side is
\begin{eqnarray*}
  &   &
\sum_{\uuu}
   |B(r,\uuu)|
\cdot (1-\yyy)^{-1} \cdot  \uuu_{1}
\cdot (1-\yyy)^{-1} \cdot  \uuu_{2}
                    \cdots \uuu_{m}
\cdot (1-\yyy)^{-1} \\
  & = &
(1-\yyy)^{-1}\cdot B\left(r;\ \xxx(1-\yyy)^{-1},
  \ \yyy(1-\yyy)^{-1}\right).
\end{eqnarray*}
\end{proof}

Combining the above results, we obtain our desired theorem:
\begin{theorem}
The generating function of the Ehrhart polynomials
of the descent polytopes is given by
$$
I(r; \xxx,\yyy) 
   =
(1-\yyy)^{-1}
   \cdot
Q_{r}\Big((\xxx-\yyy)(1-\yyy)^{-1}\Big)
   \cdot
\left(
1 
  -
\yyy(1-\yyy)^{-1}
      \cdot
Q_{r}\Big((\xxx-\yyy)(1-\yyy)^{-1}\Big)
\right)^{-1}   ,
$$
where $Q_{r}(x) = x^{-1} \cdot \left((1-x)^{-r-1} - 1\right)$.
\label{theorem_Ehrhart}
\end{theorem}

\section{Concluding remarks}

A more general invariant of the descent polytopes to study
is the flag $f$-vector. The flag $f$-vector is efficiently
encoded by the ${\bf c}{\bf d}$-index. Is there a way to describe
the ${\bf c}{\bf d}$-index of the descent polytope $\DP_S$ in terms of the
$\xxx\yyy$-word $\vvv_S$?  Finding a non-commutative generating
function for the ${\bf c}{\bf d}$-indices of descent polytopes would
be a natural way to extend the results of this paper.
The ${\bf c}{\bf d}$-indexes
of descent polytopes up to dimension $6$
can be found
in~\cite[Appendix~A.2]{Chebikin}.

Setting $t = 1$ in the polynomial $F_{\vvv}(t)$
we obtain the
number of faces of the descent polytope~$\DP_{\vvv}$.
In particular, for the alternating word $\zzz_{n}$
we obtain the sequence
$\{F_{\zzz_{n-1}}(1)\}_{n \geq 1} = 3,7,19,51,\ldots$.
This sequence has a different combinatorial interpretation,
as it matches
the sequence A052948 in the Online Encyclopedia of Integer Sequences
\cite{OEIS} defined as the number of paths from $(0,0)$ to $(n+1,0)$ with
allowed steps $(1,1)$, $(1,0)$ and $(1,-1)$ contained within the region
$-2 \leq y \leq 2$. 
The generating function
$$
\frac{1-2x^{2}}{1-3x + 2x^3}
$$
given in~\cite{OEIS} indeed
results if $t=1$ is substituted into~(\ref{equation_alternating})
and the constant $1$ is added.
Is there a bijective proof?
A first step to find such a bijective proof would
be to find a statistic on these lattice paths
with the same distribution as the dimensions
of the faces of the descent polytope
$\DP_{\zzz_{n-1}}$.

For an $\xxx\yyy$-word $\vvv = \vvv_{1} \vvv_{2} \cdots \vvv_{n}$
let $\vvv^{*}$ denote the reverse of the word, that is,
$\vvv^{*} = \vvv_{n} \cdots \vvv_{2} \vvv_{1}$.
Note that the two descent polytopes
$\DP_{\vvv}$ and $\DP_{\vvv^{*}}$
only differ by a linear transformation
and hence their $f$-polynomials agree,
that is,
$F_{\vvv} = F_{\vvv^{*}}$.
However the expressions for the
$f$-polynomials for
$F_{\vvv}$ and $F_{\vvv^{*}}$
in Corollary~\ref{corollary_factorization} differ.
Is there a bijection between the factorizations
of $\vvv$ and $\vvv^{*}$?
The number of factorizations of $\vvv$ is also equal the number
of alternating subwords of $\vvv$;
see Corollary~\ref{corollary_number_vertices}.
This fact also asks for a bijective proof.

A second way to encode subsets of $[n-1]$ is by compositions.
In~\cite[Chapter~3]{Chebikin}
this encoding is used to obtain more
recurrences to compute the $f$-polynomial $F_{S}$.

More inequalities for the descent statistic have been
proved in~\cite{Ehrenborg_Levin_Readdy,Ehrenborg_Mahajan}.
Can these inequalities be extended to the $f$-polynomial
$F_{\vvv}$? For instance, Ira Gessel asked the following question:
where does the maximum of the descent set statistic occur
when restricting to words $\vvv$ of length $n-1$ having exactly
$k$ runs of $\xxx$'s and $\yyy$'s.
He conjectured and it was proved in~\cite{Ehrenborg_Mahajan}
that the maximum occurs at the composition
$(r, \underbrace{r+1, \ldots, r+1}_{a}, r, \ldots,r)$ where
$r = \lfloor (n-1)/k \rfloor$
and
$a = (n-1) - r \cdot k$.
Would the $f$-polynomial be maximized at the 
same composition?

Descent polytopes
occur as a subdivision of the $n$-dimensional unit cube
in the work of Ehrenborg, Kitaev, and Perry~\cite{Ehrenborg_Kitaev_Perry}.
They are studying consecutive pattern avoidance
with analytic means. When considering descent pattern avoidance
they obtain operators on $L^{2}([0,1]^{n})$
whose eigenfunctions only depends on $x_{1}$ when restricted
to a descent polytope.

Theorem~\ref{theorem_Ehrhart} gives a rational
non-commutative generating function for the
Ehrhart polynomials of the descent polytopes.
However, we know that this generating function
is symmetric in the two variables $\xxx$ and $\yyy$.
Is there a different rational expression for
this generating function that shows this symmetry?

\section*{Acknowledgements}

We thank Margaret Readdy and two referees for their comments on
an earlier version of this paper.
The second author was partially funded by
National Science Foundation grant DMS-0902063
and he would like to thank the Institute for Advanced Study
where this paper was completed.

\newcommand{\journal}[6]{{\sc #1,} #2, {\it #3} {\bf #4} (#5), #6.}
\newcommand{\book}[4]{{\sc #1,} ``#2,'' #3, #4.}
\newcommand{\bookf}[5]{{\sc #1,} ``#2,'' #3, #4, #5.}
\newcommand{\books}[6]{{\sc #1,} ``#2,'' #3, #4, #5, #6.}
\newcommand{\collection}[6]{{\sc #1,}  #2, #3, in {\it #4}, #5, #6.}
\newcommand{\thesis}[4]{{\sc #1,} ``#2,'' Doctoral dissertation, #3, #4.}
\newcommand{\springer}[4]{{\sc #1,} ``#2,'' Lecture Notes in Math.,
                          Vol.\ #3, Springer-Verlag, Berlin, #4.}
\newcommand{\preprint}[3]{{\sc #1,} #2, preprint #3.}
\newcommand{\preparation}[2]{{\sc #1,} #2, in preparation.}
\newcommand{\appear}[3]{{\sc #1,} #2, to appear in {\it #3}}
\newcommand{\submitted}[3]{{\sc #1,} #2, submitted to {\it #3}}
\newcommand{\JCTA}{J.\ Combin.\ Theory Ser.\ A}
\newcommand{\AdvancesinMathematics}{Adv.\ Math.}
\newcommand{\JournalofAlgebraicCombinatorics}{J.\ Algebraic Combin.}

\bigskip
\noindent
{\em D.\ Chebikin,
     Department of Mathematics,
     MIT,
     Cambridge, MA~02139}, \\
{\tt chebikin@gmail.com}.  \\
{\em R.\ Ehrenborg,
     Department of Mathematics,
     University of Kentucky,
     Lexington, KY~40506}, \\
{\tt jrge@ms.uky.edu}.

\end{document}